\documentstyle[12pt]{article}
\input amssym.def
\input amssym.tex

\def\kk{\underline{k}}
\def\1{\underline{1}}

\def\L{\Bbb L}
\def\Z{\Bbb Z}
\def\Q{\Bbb Q}
\def\C{\Bbb C}

\newtheorem{theorem}{Theorem}
\newtheorem{statement}{Statement}

\newenvironment{definition}
{\smallskip\noindent{\bf Definition\/}:}{\smallskip\par}

\newenvironment{examples}
{\smallskip\noindent{\bf Examples\/}.}{\smallskip\par}
\newenvironment{remark}
{\smallskip\noindent{\bf Remark\/}.}{\smallskip\par}
\newenvironment{remarks}
{\smallskip\noindent{\bf Remarks\/}.}{\smallskip\par}
\newenvironment{proof}
{\noindent{\bf Proof\/}.}{{ $\square$}\smallskip\par}

\title{
An exponential function on the set of varieties
\footnote{Math. Subject Class. 14A99, 14G10}
}
\author{S.M.~Gusein-Zade \thanks{Partially supported by the grants
RFBR--01--01--00739,
INTAS--00--259, NWO--RFBR--047.008.005.
Address: Moscow State University,
Faculty of Mathematics and Mechanics, Moscow, 119899, Russia.
E-mail: sabir\symbol{'100}mccme.ru} \and
I.~Luengo
\and A.~Melle--Hern\'andez \thanks{The last two authors were partially
supported by the grant BFM2001--1488--C02--01. Address:
University Complutense de Madrid, Dept. of Algebra,
Madrid, 28040, Spain.
E-mail: iluengo\symbol{'100}mat.ucm.es, amelle\symbol{'100}mat.ucm.es}}
\date{}
\begin{document}
\def\eps{\varepsilon}

\maketitle

\begin{abstract}
Let $\cal R$ be either the Grothendieck semiring (semigroup with
multiplication) of complex quasi--projective varieties, or the Grothendieck
ring of these varieties, or the Grothendieck ring localized by the class
$\L$ of the complex affine line. We introduce a construction which defines
operations of taking powers of series over these (semi)rings. This means
that, for a power series $A(t)=1+\sum\limits_{i=1}^\infty A_i t^i$ with
the coefficients $A_i$ from $\cal R$ and for $M\in {\cal R}$, there is
defined a series $\left(A(t)\right)^M$ (also with coefficients from
$\cal R$) so that all the usual properties of the exponential function hold.
We also express in these terms the generating function of the Hilbert scheme
of points (0-dimensional subschemes) on a surface.
\end{abstract}

Let $S_0({\rm{Var}}_{\C} )$, $K_0({\rm{Var}}_{\C})$,
and $K_0({\rm{Var}}_{\C})_{(\L)}$ be the
Grothendieck semiring (semigroup with multiplication) of complex
quasi--projective varieties, the Gro\-then\-dieck ring of these varieties,
and the Grothendieck ring localized by the class $\L$ of the complex
affine line $\C$ respectively. (If we want to discuss any of them
without specification, we shall denote it by $\cal R$. In all these cases,
if $Y$ is a Zariski closed subset of $X$, then $[X]=[Y]+[X\setminus Y]$.)
Power series $\sum\limits_{i=0}^\infty A_i t^i$ with coefficients from
one of these (semi)rings are usual objects of study, in particular, in
the framework of the theory of motivic integration (see, e.g., \cite{DL},
\cite{L}). There are defined the sum and the product of such series
and therefore the $m$-th power of a series for $m\in\Z_{\ge0}$. In what follows
we shall consider power series with coefficients from one of these
(semi)rings with $A_0=1$ ($=pt$). We shall show that there exists
a (natural in some sense) notion of the $M$-th power of such a series for
the exponent $M$ from the same ring or semiring.

The problem arose from attempts to understand motivic versions of some
notions and/or statements in which the Euler characteristic participates
as an exponent of a series.
The most straightforward approach using the Taylor decomposition of the
series $(1+a_1 t+a_2 t^2+\ldots)^m$ leads to a series coefficients
of which are polynomials in $m$ and $a_i$ ($i=1, 2, \ldots$) with rational
coefficients, what has (or can be given) sense in ${\cal R}\otimes\Q$
but not in $\cal R$ itself.

The idea started from the following observation. For a complex algebraic
variety $M$, its zeta function $\zeta_M(t)$ is defined as the power series
$$
\zeta_M(t) = 1 + S^1M\cdot t + S^2M\cdot t^2 + S^3M\cdot t^3 + \ldots,
$$
where $S^kM$ is the $k$th symmetric power $M^k/S_k$ of the variety $M$,
$S_k$ is the symmetric group of  permutations on $k$ elements
(see \cite{K},\cite{L}). One has:
$\zeta_0(t) = 1$,
$\zeta_1(t) = 1+t+t^2+\ldots$,
$\zeta_{M+N}(t) = \zeta_M(t)\cdot\zeta_N(t)$
(here $0$ and $1$ are the classes of
the empty set and of a one-point set respectively). Therefore
one can consider $\zeta_M(t)$ as the exponential function $(1+t+t^2+\ldots)^M$.
(The last property above guarantees that this definition can be extended to $M$
from the Grothendieck ring $K_0({\rm{Var}}_{\C})$:
$\zeta_{-M}(t):=\left(\zeta_M(t)\right)^{-1}$ (of course in this case the
coefficients of the series have to be considered as elements of the ring $K_0({\rm{Var}}_{\C})$ as
well).) One can say that the aim of this paper is to extend this definition of a
power of the series $1+t+t^2+\ldots$ to an arbitrary series
$1+A_1 t+A_2 t^2+\ldots$ with the coefficients $A_i$ from
one of the (semi)rings under discussion ($S_0({\rm{Var}}_{\C})$, $K_0({\rm{Var}}_{\C})$, or
$K_0({\rm{Var}}_{\C})_{(\L)}$).

Let $1+{\cal R}_+[[t]]$ be the set of (formal) power series $A(t)=1+\sum\limits_{i=1}^\infty A_it^i$
with $A_i\in{\cal R}$. If $\cal R$ is a ring, $1+{\cal R}_+[[t]]$ is an abelian group
with respect to the multiplication. If $\cal R$ is a semiring, $1+{\cal R}_+[[t]]$ is an (abelian)
semigroup (in this case the division is not defined).

\begin{definition}
An operation of taking powers of series over a (semi)ring $\cal R$ is a map from
$\left(1+{\cal R}_+[[t]]\right)\times {\cal R}$  to $1+{\cal R}_+[[t]]$ which is denoted by
$\left(A(t)\right)^M$ ($A(t)= 1+\sum\limits_{i=1}^\infty A_it^i \in 1+{\cal R}_+[[t]]$,
$A_i\in{\cal R}$, $M\in{\cal R}$) and which possesses the properties:
\begin{enumerate}
\item $\left(A(t)\right)^0=1$.
\item $\left(A(t)\right)^1=A(t)$.
\item $\left(A(t)\cdot B(t)\right)^M=\left(A(t)\right)^M\cdot \left(B(t)\right)^M$.
\item $\left(A(t)\right)^{M+N}=\left(A(t)\right)^M\cdot \left(A(t)\right)^N$.
\item $\left(A(t)\right)^{MN}=\left(\left(A(t)\right)^N\right)^M$.
\end{enumerate}
\end{definition}

\begin{remark}
In what follows we shall permit ourselves to denote $\left(A(t)\right)^M$
by $A^M(t)$ as well. This is convenient in cases when we substitute the variable $t$
by an expression (e.g., by $t^s$ or by $\L t$ below). For instance, the meaning of the expression $A^M(t^s)$ is clear
and is different from the (natural) meaning of $\left(A(t^s)\right)^M$. (In fact
for our construction these two expressions happen to coincide, however, this follows
only from the definition bellow.)
Usually we shall
use the same notations for a variety (or for a constructible set) and for its class
in the corresponding semiring or ring.
\end{remark}

\begin{theorem}\label{theo1}
There exists an operation of taking powers of series over the Grothendieck semiring
$S_0({\rm{Var}}_{\C})$ of complex algebraic varieties such that
$$
(1+t+t^2+\ldots)^M=1+S^1M\cdot t+S^2M\cdot t^2+\ldots=\zeta_M(t).
$$
\end{theorem}

\begin{proof}
Before giving the construction (definition) of the operation, we shall try to explain where it comes from.
One property which we would like to have is
$$
\left(A(t)\right)^{MN} =
\left(\left(A(t)\right)^N\right)^M.
$$
In particular, this means that
$$
1 + S^1(MN) t + S^2(MN) t^2 + \ldots=\zeta_{MN}(t)=(1+t+t^2+\ldots)^{MN}=
$$
$$
=\left(\zeta_N(t)\right)^M=
(1 + S^1N\cdot t + S^2N\cdot t^2 + \ldots)^M.
$$
The symmetric power $ S^k(MN)$ of the product $MN$ can be obtained from the
symmetric powers $S^iN$ and the variety $M$ in the following way. Let
$\pi:(MN)^k\to M^k$ and $\widetilde\pi:S^k(MN)\to S^kM$ be the natural
projections. The spaces $ M^k$ and $S^kM$ have natural decompositions
into parts corresponding to different patterns of coinciding points.
For $\kk=(k_1, k_2, \ldots)$ with $\sum\limits_i ik_i=k$, let $M_{\kk}$
(respectively $\widetilde M_{\kk}$) be the subset of $ M^k$ (respectively
of $S^kM$), points of which ($k$-tuples of points of $M$, ordered or unordered
respectively) consist of $k_1$ different points of $M$, $k_2$ (different) pairs
of equal points, $k_3$ triples of equal points, {\dots} The subspace $M_{\kk}$
is locally closed in $M^k$, $M^k=\bigcup\limits_{\kk:\Sigma ik_i=k}M_{\kk}$,
$\widetilde M_{\kk}=M_{\kk}/S_k$,
$S^kM=\bigcup\limits_{\kk:\Sigma ik_i=k}\widetilde M_{\kk}$. The preimage
$\widetilde\pi^{-1}(\widetilde M_{\kk})$ of the set $\widetilde M_{\kk}$
in $S^k(MN)$ can be obtained in the following way. Consider the direct product
$$
V_{\kk}=M_{\kk}\times {\underbrace{N\times\ldots\times N}_{k_1}} \times
{\underbrace{S^2N\times\ldots\times S^2N}_{k_2}}\times\ldots
$$
One has $M_{\kk}\cong \left(\prod\limits_i M^{k_i}\right)\setminus \Delta
\,\left(=M^{\Sigma k_i}\setminus\Delta\right)$, where $\Delta$ is the
"large diagonal" in $M^{\Sigma k_i}$ which consists of $(\sum k_i)$-tuples
of points of $M$ with at least two coinciding ones. There is a natural
free action of the group $\prod\limits_i S_{k_i}$ on the space $V_{\kk}$,
where $S_{k_i}$ acts by permuting corresponding $k_i$ factors in
$\prod\limits_i M^{k_i}\supset M_{\kk}$ and the spaces $S^iN$
(there are $k_i$ of them in the product). The
preimage $\widetilde\pi^{-1}(\widetilde M_{\kk})$ is nothing else but
the factor of the space $V_{\kk}$ by this action, $S^k(MN)=
\bigcup\limits_{\kk:\Sigma ik_i=k}\widetilde\pi^{-1}(\widetilde M_{\kk})$.
To get the definition of the coefficients of the series
$\left(A(t)\right)^M=(1+A_1 t+A_2 t^2+\ldots)^{M}$, one should substitute the spaces
$S^iN$ in the description above by the spaces $A_i$. This leads to the
following definition.

\begin{definition}
$$
\left(A(t)\right)^M=(1+A_1 t+A_2 t^2+\ldots)^M=
$$
$$
=1+\sum_{k=1}^\infty
\left
\{\sum_{\kk:\sum ik_i=k}
\left(
\left(
(\prod_i M^{k_i}
)
\setminus\Delta
\right)
\times\prod_i A_i^{k_i}/\prod_iS_{k_i}
\right)
\right\}
\cdot t^k,
$$
where the group $S_{k_i}$ acts by permuting corresponding $k_i$ factors in
$\prod\limits_i M^{k_i}\supset (\prod_i M^{k_i})\setminus\Delta$ and the spaces $A_i$.
\end{definition}

\begin{remark}
This definition can also be "read from" the following formula for the (usual) power
$(1+ a_1 t+ a_2 t^2+ \ldots)^m$ of a series with a natural exponent:
$$
(1+a_1 t+a_2 t^2+\ldots)^m=
$$
$$
=1+\sum_{k=1}^\infty
\left(
\sum_{\kk:\Sigma ik_i=k}
m(m-1)\ldots(m-\Sigma k_i +1)
\cdot\prod_i a_i^{k_i}/\prod_i k_i!
\right)
\cdot t^k
$$
(see, e.g., \cite{St}, page 40).
In this formula one should understand the product $m(m-1)\ldots(m-\Sigma k_i +1)$
as $M^{\Sigma k_i}\setminus\Delta$ (this product is the number of elements of the
indicated space if M is a set which consists of $m$ points) and should understand the division
by $\prod\limits_i k_i!$ as the factorization by the group $\prod\limits_iS_{k_i}$.
\end{remark}

The fact that $(1+t+t^2+\ldots)^M=1+S^1M\cdot t+S^2M\cdot t^2+\ldots=\zeta_M(t)$
and the properties 1 and 2 from the definition of an operation of taking powers of series
are obvious. Let us reformulate the definition of the series
$\left(A(t)\right)^M$ a little bit so that the properties 3 to 5 will
be proved by establishing one-to-one correspondences between sets representing the
coefficients of the left hand side ({\bf LHS}) and the right hand side ({\bf RHS})
series. (We don't use this version as the initial definition since in this case one
has to define the corresponding sets as constructible ones.)
The coefficient at the monomial $t^k$
in the series $\left(A(t)\right)^M$ is represented by the set whose element
is a (finite) subset $K$ of points of the variety $M$ with (positive)
multiplicities such that
the total number of points of the set $K$ counted with multiplicities is equal to $k$
plus a map $\varphi$ from $K$ to ${\cal A}=\coprod\limits_{i=0}^\infty A_i$ ($A_0=1=pt$) such that
a point of multiplicity $s$ goes to $A_s\subset{\cal A}$ ($\coprod$ means the disjoint union,
i.e., the sum in the semiring $S_0({\rm{Var}}_{\C})$).
For short, instead of writing that "a coefficient of a series is represented by a set which
consists of elements of the form\dots" we shall write that "an element of the coefficient is\dots".

\underline{Proof of 3.} Let $C_j$ be the coefficient at the monomial $t^j$ in the product
$A(t)\cdot B(t)$ ($B(t)=\sum\limits_{i=0}^\infty B_it^i$): $C_j=\sum\limits_{i=0}^j A_iB_{j-i}$.
An element of the
coefficient at the monomial $t^k$
in the {\bf LHS} of the equation 3 is a $k$-point subset $K$ of $M$
with a map $\varphi$ from $K$ to $\coprod\limits_{i=0}^\infty C_i=
\coprod\limits_{i,j\ge 0} A_iB_j$ such that a point of multiplicity $s$ goes to
$C_s=\coprod\limits_{i=0}^s A_i\times B_{s-i}$. An element of the coefficient at the monomial $t^k$ in the {\bf RHS}
of this equation is: an $\ell$-point subset $K_1$ of the variety $M$ with $0\le\ell\le k$,
a map $\varphi_1$ from $K_1$ to ${\cal A}=\coprod\limits_{i=0}^\infty A_i$ (such that
a point of multiplicity $s$ goes to $A_s$), a $(k-\ell)$-point subset $K_2$ of the variety $M$,
a map $\varphi_2$ from $K_2$ to ${\cal B}=\coprod\limits_{i=0}^\infty B_i$ (such that
a point of multiplicity $s$ goes to $B_s$). Suppose we have an element of the first set.
Let us decompose the subset $K$ into two parts $K_1$ and $K_2$. A point $x$ of $K$ of multiplicity
$s\ (=s(x))$ goes to $C_s=\coprod\limits_{i=0}^s A_i\times B_{s-i}$ and thus $\varphi(x)$
belongs to one of the summands, say to
$A_{i_0}\times B_{s-i_0}$. Let us include the point $x$ into the set $K_1$ with the multiplicity
$i_0$ and into the set $K_2$ with the multiplicity $s-i_0$. If $i_0$ and/or $s-i_0$ is positive,
define $\varphi_1(x)$ and/or $\varphi_2(x)$ as the corresponding projection (to $A_{i_0}$ and/or
to $B_{s-i_0}$ respectively) of the point $\varphi(x)\in A_{i_0}\times B_{s-i_0}$. In the other direction, from an
element of the second set one can construct an element of the first one uniting the subsets
$K_1$ and $K_2$ (so that a multiplicity of a point $x$ in $K_1\cup K_2$ is equal to the sum $s_1+s_2$
of its multiplicities in $K_1$ and in $K_2$) and defining $\varphi(x)$ as $\left(\varphi_1(x),\varphi_2(x)\right)\in
A_{s_1}\times B_{s_2}\subset C_{s_1+s_2}$. One easily sees that these correspondences are
inverse to each other and thus are one-to-one.

\underline{Proof of 4.} An element of the coefficient at the monomial $t^k$
in the {\bf LHS} of the equation is a $k$-point subset $K$ of the variety $M\cup N$
with a map $\varphi$ from $K$ to ${\cal A}=\coprod\limits_{i=0}^\infty A_i$.
An element of the coefficient at the monomial $t^k$ in the {\bf RHS} is a pair:
an $\ell$-point subset $K_1$ of the variety $M$ with $0\le\ell\le k$ and
a $(k-\ell)$-point subset $K_2$ of the variety $N$ with their maps to ${\cal A}$.
The union of these two subsets with the corresponding map to ${\cal A}$
gives a subset in the union $M\cup N$ with a map to ${\cal A}$. This
correspondence is obviously one-to-one.

\underline{Proof of 5.} An element of the coefficient at the monomial $t^k$
in the {\bf LHS} of the equation is a $k$-point subset $K$ of the product $M\times N$
with a map $\varphi:K\to {\cal A}=\coprod\limits_{i=0}^\infty A_i$.
An element of the coefficient at the monomial $t^k$ in the {\bf RHS} is a $k$-point subset $K_M$ of
the variety $M$ with a map from it which sends a point of multiplicity $s$ to an
$s$ point subset of $N$ with a map from it to ${\cal A}$. To establish a one-to-one correspondence
between these coefficients (or rather between the variety representing them)
one should define $K_M$ as the projection $\mbox{pr}_M K$ of the subset $K\subset M\times N$
to the first factor and, for $x\in K_M$, the corresponding subset in $N$ as
$\mbox{pr}_N \mbox{pr}^{-1}_M(x)$ (with the natural map to ${\cal A}$).
\end{proof}

\begin{examples}
$$
(1 + t + t^2 + ...)^M = 1 + \sum_{k=1}^\infty S^kM\cdot t^k;
$$
$$
(1 + t)^M = 1 + \sum_{k=1}^\infty (M^k\setminus\Delta)/S_k\cdot t^k.
$$
\end{examples}

Let $\chi(X)$ be the Euler characteristic of the space $X$. For a series $A(t)=1+ A_1 t+ A_2t^2+\ldots$,
one defines its Euler characteristic as the series
$$
\chi\left(A(t)\right) =1+ \chi(A_1)t+ \chi(A_2)t^2+\ldots\,\in\, \Z[[t]].
$$

\begin{statement}\label{st1}
$$
\chi\left(\left(A(t)\right)^M\right)=\left(\chi\left(A(t)\right)\right)^{\chi(M)}.
$$
\end{statement}

The {\bf proof} follows either from direct calculations or from the fact that
the coefficients at the monomials $t^{k}$ in the {\bf LHS} of the equation
are polynomials in the Euler characteristics of the varieties
$M$ and $A_i$ ($i=1,2,\ldots)$ and the equation holds for "natural numbers", i.e.,
 for the case when all the varieties $M$ and $A_i$ are finite sets of points. { $\square$}\smallskip

\begin{remark}
The corresponding formula for the virtual Hodge--Deligne polynomial
is not so straightforward.
\end{remark}

\begin{theorem}\label{theo2}
There exists an operation of taking powers of series over the Grothendieck ring
$K_0({\rm{Var}}_{\C})$ of complex algebraic varieties which extends the one
defined over the semiring $S_0({\rm{Var}}_{\C})$.
\end{theorem}

\begin{proof}
To define the operation notice that for any series $A(t)\in 1+ K_0({\rm{Var}}_{\C})_+[[t]]$
there exists a series $B(t)=\sum\limits_{i=0}^\infty B_it^i\in 1+ K_0({\rm{Var}}_{\C})_+[[t]]$ with
the coefficients from the image of the natural map $S_0({\rm{Var}}_{\C})\to K_0({\rm{Var}}_{\C})$ such that
all the coefficients of the product $C(t)=A(t)\cdot B(t)$ are from the same image as well.
Then one puts:
$$
\left(A(t)\right)^M:=\left(C(t)\right)^M/\left(B(t)\right)^M.
$$
To define the power of a series with the exponent $M$ from the Grothendieck ring $ K_0({\rm{Var}}_{\C})$,
one puts
$$
\left(A(t)\right)^{-M}:=1/\left(A(t)\right)^M.
$$
The properties of the definition of an operation of taking powers of series obviously hold.
\end{proof}

\begin{theorem}\label{theo3}
There exists an operation of taking powers of series over the ring $K_0({\rm{Var}}_{\C})_{(\L)}$
$($i.e., over the Grothendieck ring of complex algebraic varieties localized by the class $\L$
of the complex affine line $\C$$)$ which extends the one
defined over the ring $K_0({\rm{Var}}_{\C})$.
\end{theorem}

\begin{proof}
First let us define the operation for series $A(t)$ with the coefficients $A_i$ from the ring $K_0({\rm{Var}}_{\C})_{(\L)}$
with the exponent $M$ from the "non-localized" ring $K_0({\rm{Var}}_{\C})$.
This is possible because of the following statement.

\begin{statement}\label{st2}
Let $A_i$ and $M$ be from the Grothendieck ring of complex algebraic varieties. Then, for $s\ge0$,
$\left(A(\L^s t)\right)^M = A^M(\L^s t)$.
\end{statement}

\begin{proof}
It is sufficient to prove this equation for a series with the coefficients from $S_0({\rm{Var}}_{\C})$.
The coefficient at the monomial $t^k$ in the power series $A^M(t)$ is a sum of varieties of the form
$$
V=\left(
(\prod_i M^{k_i}
)
\setminus\Delta
\right)
\times\prod_i A_i^{k_i}/\prod_iS_{k_i}
$$
with $\sum ik_i=k.$ The corresponding summand $\widetilde V$ in the coefficient
at the monomial $t^k$ in the power series $\left(A(\L^s t)\right)^M$
has the form
$$\widetilde V=\left(
(\prod_i M^{k_i}
)
\setminus\Delta
\right)
\times\prod_i (\L^{si}A_i)^{k_i}/\prod_iS_{k_i}.
$$
There is a natural map $\widetilde V \to V$
which from the point of view of differential geometry is a complex analytic
vector bundle of rank $sk$ (i.e., it is locally trivial over a neighbourhood of each point in the usual topology). According to \cite{S}
this is a vector bundle in the "algebraic sense", i.e., it is locally trivial over a Zariski open neighbourhood of each point. This implies that $[\widetilde V]=\L^{sk}\cdot[V].$
\end{proof}

For a series $A(t)=\sum\limits_{i=0}^\infty A_it^i$, let $J^r A(t)$ be its truncation up to terms of degree $r$, i.e.,
$J^r A(t)=\sum\limits_{i=0}^r A_it^i$.
Statement \ref{st2} implies that, for $A_i\in K_0({\rm{Var}}_{\C})_{(\L)},$
$i=1,2,\ldots$, $M\in K_0({\rm{Var}}_{\C})$, one can define the power
$\left(A(t)\right)^M$ by the formula
$$
J^r\left(\left(A(t)\right)^M\right):=
\left(J^r A(\L^s t)\right)^M
\mbox{\raisebox{-0.5ex}{$\vert$}}{}_{t\mapsto {t}/{\L^s}}
$$
for $s$ large enough so that all the coefficients of $J^r A(\L^s t)$ belong to the
image of the map $ K_0({\rm{Var}}_{\C})\to  K_0({\rm{Var}}_{\C})_{(\L)}.$
One can easily see that the properties 3--5 hold.

Now we have to extend the operation to the exponent $M$ from the localized ring
$K_0({\rm{Var}}_{\C})_{(\L)}.$ First let us do it for one particular series, namely for $(1+t+t^2+\ldots).$ One can say that we define the zeta function $\zeta_M(t)$
for $M \in K_0({\rm{Var}}_{\C})_{(\L)}.$

\begin{statement}\label{st3}
For  $M \in K_0({\rm{Var}}_{\C})$, $s\ge 0$,
$$
\zeta_{\L^s M}(t)=\zeta_{M}(\L^s t).
$$
\end{statement}

\begin{proof}
It is sufficient to prove the equation
for $s=1$. One has $\zeta_{\L}(t)=1+\L t+\L^2 t^2+\ldots$ and therefore
$$
\zeta_{\L M}(t)=(1+\L t+\L^2 t^2+\ldots)^M=(1+t+ t^2+\ldots)^M
\mbox{\raisebox{-0.5ex}{$\vert$}}{}_{t\mapsto {\L t}},
$$
(we have used Statement~\ref{st2}).
\end{proof}

This statement permits to define $\zeta_M(t)$ for $M \in K_0({\rm{Var}}_{\C})_{(\L)}$
by the formula
$$
\zeta_{M}(t)=\zeta_{\L^s M}(\L^{-s} t)
$$
for $s$ large enough (so that $\L^s M$ belongs to the image of the map
$K_0({\rm{Var}}_{\C})\to K_0({\rm{Var}}_{\C})_{(\L)}$).

For $A(t)\in 1+ K_0({\rm{Var}}_{\C})_+[[t]]$, $M\in K_0({\rm{Var}}_{\C})$, $s>0$, one has
$$
\left(A(t^s)\right)^M=A^M(t^s).
$$
Since $\zeta_M(t)=1+Mt+\ldots$, for any series $A(t)\in 1+ K_0({\rm{Var}}_{\C})_{(\L)+}[[t]]$
and for any $r>0$, the truncated series $J^rA(t)$ can be (in a unique way) represented as
$$
J^r\left(\zeta_{A_1}(t)\cdot\zeta_{A_2'}(t^{2})\cdot\zeta_{A_3'}(t^{3})\cdot\ldots\cdot
\zeta_{A_r'}(t^{r})\right)
$$
with $A_i'\in  K_0({\rm{Var}}_{\C})_{(\L)}.$
One defines
$\left(A(t)\right)^M$ by the formula
\begin{equation}\label{formul}
J^r\left(\left(A(t)\right)^M\right):=J^r\left(\zeta_{MA_1}(t)\cdot\zeta_{MA_2'}(t^{2})\cdot\ldots\cdot
\zeta_{MA_r'}(t^{r})\right).
\end{equation}
Properties 3--5 of the definition of taking powers of series obviously hold.
\end{proof}

\begin{remarks}
{\bf 1.} If the zeta function $\zeta_M(t)$ itself is already defined for
$M$ from $K_0({\rm{Var}}_{\C})$ or
from $K_0({\rm{Var}}_{\C})_{(\L)}$, the formula (\ref{formul})
can be used to define the operation of taking powers of series over
these rings. Such a operation is {\bf unique} (because of Property 6).
However it does not help to define this operation over $S_0({\rm{Var}}_{\C})$.
\newline
{\bf 2.} The described structure is somewhat similar to the
structure of a $\lambda$-ring (see, e.g., \cite{Knut}), but is different
from it. If the series $\lambda_t(M)=(1+t)^M$ would define the structure of a
$\lambda$-ring, then the coefficient at $t^2$ in $(1+t)^{\L^2}$
(computed as $\lambda_t(\L\times\L)$) has to be equal to $2\L^4-2\L^2$
when in our construction it is equal to $\L^4-\L^2$.
\newline
{\bf 3.} Yu.I.Manin has explained us that his calculations in
\cite{M} can be considered as a particular case of our construction.
\end{remarks}

\medskip

{\bf An application.} For a smooth quasi--projective surface $M$,
let $M^{[n]}={\rm{Hilb}}^nM$ be
the Hilbert scheme of $0$-dimensional subschemes of length $n$ on $M$.

\begin{statement}
In the Grothendieck ring $K_0({\rm{Var}}_{\C})$ one has:
\begin{eqnarray*}
1+\sum_{n\geq 1} M^{[n]} t^n & = & \prod\limits_{k\geq 1}
\zeta_{\L^{-1}M} ((\L t)^k)\\
& =& \left(\prod\limits_{k\geq 1} \zeta_{M} (\L^{k-1}t^k)\right)\\
& =& \left(\prod\limits_{k\geq 1} \frac{1}{1-\L^{k-1}t^k}\right)^{M}\\
& =& \prod\limits_{k\geq 1} \left(\frac{1}{1-t^k}\right)^{\L^{k-1}M}.
\end{eqnarray*}
\end{statement}

The {\bf proof} follows from the following result of L.~G\"ottsche
\cite{G}:
$$
M^{[n]} = \sum_{\kk:\Sigma ik_i=n}S^{\kk} M\cdot \L^{n-\vert\kk\vert},
$$
where $\vert\kk\vert=\sum k_i$, $S^{\kk} M=
S^{k_1}M\times\cdots\times S^{k_n}M$. Therefore
\begin{eqnarray*}
& 1+\sum\limits_{n\geq 1} M^{[n]} t^n =
\,\sum\limits_{n\geq 0} \left(\sum\limits_{\kk:\Sigma ik_i=n}
S^{\kk}M\cdot\L^{n-\vert\kk\vert}\right)\cdot t^n \\
=& \,\sum\limits_{n\geq 0}
 \left(\sum\limits_{\kk:\Sigma ik_i=n}S^{k_1}M\cdot\ldots\cdot S^{k_n}M\cdot
 \L^{n-(k_1+\ldots+k_n)}\right)\cdot t^{\Sigma ik_i}\\
=& \,\sum\limits_{n\geq 0} \left(\sum\limits_{\kk:\Sigma ik_i=n}
(S^{k_1}M\cdot \L^{-{k_1}})\cdot\ldots\cdot (S^{k_n}M\cdot
 \L^{-{k_n}})\right)\cdot(\L t)^{\Sigma ik_i} \\
=& \, \prod\limits_{k\geq 1}\left(\sum\limits_{r\ge 0} S^rM\cdot{\L^{-r}}
 (\L t)^{kr}\right)=\prod\limits_{k\geq 1} \zeta_{\L^{-1}M} \left((\L t)^k\right).\quad\square
\end{eqnarray*}

\begin{remark}
D.~van Straten told us that
the fact that the generating function of the Hilbert scheme
of 0-dimensional subschemes of a surface is in some sense an
exponent was conjectured by several people.
\end{remark}

\smallskip
{\bf Acknowledgements:} The authors are thankful to Tom\'as L. G\'omez
for useful discussions.

\end{document}